\documentclass[11pt]{amsart2000}
\newtheorem*{theorem}{Theorem}
\newtheorem*{corollary}{Corollary}
\theoremstyle{definition}
\newtheorem*{remark}{Remark}

\begin{document}
\setlength{\unitlength}{0.01in}
\linethickness{0.01in}
\begin{center}
\begin{picture}(474,66)(0,0)
\multiput(0,66)(1,0){40}{\line(0,-1){24}}
\multiput(43,65)(1,-1){24}{\line(0,-1){40}}
\multiput(1,39)(1,-1){40}{\line(1,0){24}}
\multiput(70,2)(1,1){24}{\line(0,1){40}}
\multiput(72,0)(1,1){24}{\line(1,0){40}}
\multiput(97,66)(1,0){40}{\line(0,-1){40}}
\put(143,66){\makebox(0,0)[tl]{\footnotesize Proceedings of the Ninth Prague Topological Symposium}}
\put(143,50){\makebox(0,0)[tl]{\footnotesize Contributed papers from the symposium held in}}
\put(143,34){\makebox(0,0)[tl]{\footnotesize Prague, Czech Republic, August 19--25, 2001}}
\end{picture}
\end{center}
\vspace{0.25in}
\setcounter{page}{91}
\title{Transfinite sequences of continuous and Baire~1 functions on 
separable metric spaces}
\author{M\'arton Elekes}
\address{Department of Analysis, E\"otv\"os Lor\'and University\\
Budapest, P\'azm\'any P\'eter s\'et\'any 1/c, 1117, Hungary}
\email{emarci@cs.elte.hu}
\thanks{This is a research announcement. A complete article, written
jointly with Kenneth Kunen, will be published elsewhere}
\thanks{M\'arton Elekes,
{\em Transfinite sequences of continuous and Baire 1 functions on 
separable metric spaces},
Proceedings of the Ninth Prague Topological Symposium, (Prague, 2001),
pp.~91--92, Topology Atlas, Toronto, 2002}
\begin{abstract}
We investigate the existence of well-ordered sequences of Baire 1 
functions on separable metric spaces.
\end{abstract}
\maketitle

Any set $\mathcal{F}$ of real valued functions defined on an arbitrary
set $X$ is partially ordered by the pointwise order, that is $f\leq g$ iff
$f(x)\leq g(x)$ for all $x\in X$. In other words put $f<g$ iff $f(x)\leq
g(x)$ for all $x\in X$ and $f(x)\not= g(x)$ for at least one $x\in X$. 
Our aim will be to investigate the possible length of the increasing or
decreasing well-ordered sequences of functions in $\mathcal{F}$ with
respect to this order. 

A classical theorem of Kuratowski asserts, that if $\mathcal{F}$ is the
set of continuous or Baire~1 functions defined on a Polish space $X$, then
there exists a monotone sequence of length $\xi$ in $\mathcal{F}$ iff
$\xi<\omega_1$ (see \cite[\S 24. III.2']{Ku}). 
Moreover, P. Komj\'ath proved in \cite{Ko} that the corresponding question
concerning Baire~$\alpha$ functions for $2\leq\alpha<\omega_1$ is
independent of $ZFC$. 

In the present paper we investigate what happens if we drop the condition
of completeness and replace the Polish space $X$ by a separable metric space.

Our main results are the following. Let $d(X)$ denote the density of a space
$X$.

\begin{theorem}
Let $(X,\varrho)$ be a metric space. 
Then there exists a well-ordered sequence of length $\xi$ of continuous
real-valued functions defined on $X$ iff $\xi<d(X)^+$.
\end{theorem}

\begin{corollary}
A metric space is separable iff every well-ordered sequence of continuous
functions defined on it is countable.
\end{corollary}

\begin{theorem}
There exists a separable metric space on which there exists a well-ordered
sequence of length $\omega_1$ of Baire~1 functions.
\end{theorem}

\begin{theorem}
The following statement: `There exists a separable metric space on which
there exists a well-ordered sequence of length $\omega_2$ of Baire~1
functions' is independent of $ZFC+\neg CH$.
\end{theorem}

\begin{remark} 
During and after the conference Kenneth Kunen answered one of my
questions, and also improved some of the results and proofs. 
These results will appear in a forthcoming joint paper.
\end{remark}

%\bibliographystyle{amsplain}
%\bibliography{09}
\providecommand{\bysame}{\leavevmode\hbox to3em{\hrulefill}\thinspace}
\providecommand{\MR}{\relax\ifhmode\unskip\space\fi MR }
\providecommand{\MRhref}[2]{%
  \href{http://www.ams.org/mathscinet-getitem?mr=#1}{#2}
}
\providecommand{\href}[2]{#2}

\end{document}